\documentclass[12pt]{article}
\usepackage[cp1251]{inputenc}
\usepackage[T2A]{fontenc}
\usepackage[english]{babel}
\usepackage{amsfonts}
\usepackage{amsmath}
\usepackage{amssymb}
\usepackage{graphicx}
\usepackage[epsf]{}
\usepackage{amscd}
\usepackage{amsthm}

\newcommand{\R}{{\mathbb R}}

\newtheorem{theorem}{Theorem}

\newtheorem{proposition}{Proposition}
\newtheorem{lemma}{Lemma}

\newtheorem{Rm}{Remark}                                     %
\newcommand{\e}{\end{equation}}
\renewcommand{\b}{\begin{equation}}
\newcommand{\eps}{\varepsilon}

\title{On interrelations between \\divergence-free and Hamiltonian dynamics}
\author{L. Lerman, E. Yakovlev\\
\normalsize
Lobachevsky State University of Nizhny Novgorod, Russia\\
}
\date{}
\begin{document}
\maketitle

\begin{abstract}
A mathematically correct description is presented on the interrelations
between the dynamics of divergence free vector fields on an oriented
3-dimensional manifold $M$ and the dynamics of Hamiltonian systems. It is
shown that for a given divergence free vector field $X$ with a global cross-section
there exist some 4-dimensional symplectic manifold $\tilde{M}\supset M$ and
a smooth Hamilton function $H: \tilde{M}\to \R$ such that for some $c\in \R$ one gets $M = \{H=c\}$
and the Hamiltonian vector field $X_H$ restricted on this level coincides with
$X$. For divergence free vector fields with singular points such the
extension is impossible but the existence of local cross-section allows
one to reduce the dynamics to the study of symplectic diffeomorphisms in
some sub-domains of $M$. We also consider the case of a divergence free vector field
$X$ with a smooth integral having only finite number of critical levels. It is
shown that such a noncritical level is always a 2-torus and restriction of $X$ on
it possesses a smooth invariant 2-form. The linearization of the flow on
such a torus (i.e. the reduction to the constant vector field) is not
always possible in contrast to the case of an integrable Hamiltonian system
but in the analytic case ($M$ and $X$ are real analytic), due to the Kolmogorov theorem, such the
linearization is possible for tori with Diophantine rotation numbers.
\end{abstract}

\section{Introduction}

It is a rather frequent case when the Lagrangian description of liquid flows
discovers structures characteristic for Hamiltonian dynamics \cite{Lamb,Ar,MW,BSh}.
Our aim in this note is to display the reason of this in the explicit
form. In a sense, we proceed what was done by Arnold in \cite{Ar,Arn2}.
To show the interrelations between these two types of dynamics we present some
calculations which have to demonstrate how symplectic 2-dimensional maps
arise in the Lagrangian description of liquid flows.
One has to say that models, where divergence-free 3-dimensional vector
fields are studied, appear not only at the Lagrangian description of liquid currents.
Such equations arise also in the study of motions of thin liquid films flowing down
along the inclined plane \cite{BSh}. The same features are characteristic
in other situations where the divergence-free flows arise \cite{AKh}. This
is observed in models of magnetic hydrodynamics, plasma confinement problems where
magnetic line fields are studied  \cite{HM,F,Gr,Kozlov}. An important application
in plasma physics involves integration of magnetic field lines
$$
\frac{dx}{dt_0} = B(x),
$$
where $t_0$ is an artificial time-like parameter parameterizing
motion along the field line, and  $div B$ is equal to
zero \cite{HM}. In any case, the study of divergence-free vector fields is
a very interesting problem itself, many features of their flow orbit
structure deserve of a detailed investigation, see, for instance \cite{Fer,NVV}.

The Lagrangian description of stationary 3-dimensional flows in the space
$\mathbb R^3$ with coordinates $(x,y,z)$ has the form
$$
\dot x = A(x,y,z,t),\;\dot y = B(x,y,z,t),\;\dot z = C(x,y,z,t),
$$
where vector field $V=(A,B,C)$ is the velocity field of liquid particles. If the
liquid is incompressible, then $div V \equiv 0$ and the field is
divergence-free. We consider here the autonomous case and suppose also that
vector field $V$ has some (possibly local) cross-section,
that is, there is a 2-dimensional smooth submanifold $N$ such that $N$ is
transverse to the vector field at points of $N$ and for which orbits
starting at points of some its subdomain $N_1 \subset N$ return to $N$ in finite
times.

It is more convenient to carry out the calculations in the invariant (coordinate-free) form \cite{AKh}.
Consider a smooth ($C^\infty$) 3-dimensional oriented manifold $M$ with a smooth volume form $\Omega$.
Recall that for each smooth vector field $X$ on $M$ a smooth function is defined called
{\em the divergence of the vector field}, $div\,X$. This function is given by the relation
$$
L_X\Omega = d(\imath_X\Omega)=(div\,X)\Omega,
$$
where $\imath_X\Omega$, called the interior product of $\Omega$ and $X$, is the 2-form such that for any $m\in M,$ $\xi,\eta \in T_mM$
gives $\imath_X\Omega(\xi,\eta)= \Omega(X,\xi,\eta)$. If $div\,X \equiv 0$, then they say the vector field be divergence-free.
For any given diffeomorphism $f:M \to M$ the pullback form $f^*\Omega$ is defined as $[f^*\Omega]_m(\xi,\eta,\zeta)=$
$\Omega_{f(m)}(Df(\xi),Df(\eta),Df(\zeta))$. One says that $f$ preserves
the volume if the identity $f^*\Omega = \Omega$ holds true.

Let a vector field $X$ be divergence-free. Then the flow for this vector field, i.e. one-parameter group of diffeomorphisms
$f^t$ (shifts along flow orbits) consists of divergence-free diffeomorphisms \cite{AKh}.

A nondegenerate volume form $\Omega$ on an oriented manifold $(M,\Omega)$ generates the correspondence between
vector fields $X$ on $M$ and 2-forms $\omega$. This is given by the relation
\begin{equation}\label{for}
\omega_X = \imath_X\Omega,
\end{equation}
for any given vector field $X$ on $M$.

A standard example of a 3-dimensional divergence-free vector field (DFVF, for brevity) is a
flow on a nonsingular level of a Hamilton function $H$ for a Hamiltonian vector field $X_H$ on some
smooth symplectic 4-dimensional manifold $(N,\omega)$ with the Hamiltoninan
$H: N\to \R$ and symplectic 2-form $\omega$. Suppose $dH \ne 0$ on some level $H=c$, then this level is a
smooth 3-dimensional orientable submanifold $M$ of $N$ (nonsingular level). The flow in $N$ generated by the vector
field $X_H$ preserves 2-form $\omega$ and hence, the iterated volume form $\Omega = \omega\wedge
\omega$ as well. Let us endow $N$ with some smooth Riemannian metric $g$. Then a vector field $\nabla{H}$ on $M$
is correctly defined by the relation $\imath_{\nabla{H}}g = g(\nabla H,\cdot)=dH(\cdot)$.
By construction, one gets $\nabla{H}\ne 0$ on $M$. Therefore, the function
$\rho=1/g(\nabla{H},\nabla{H})$ and a vector field $n=\rho\nabla{H}$ are defined and differentiable. Observe that
\begin{equation}
dH(n)=(\imath_{\nabla{H}}g)(n)=g(\nabla{H},n)\equiv 1.
\end{equation}
Let us set $\Omega_n=(\imath_n\Omega)|_{M}$.
It turns out that $\Omega_n$ is
nondegenerate, i.e. is a volume form and invariant w.r.t. the flow on $M$.
The invariant proof of this assertion is given below in Sec. \ref{int}.
In statistical physics the measure on $M$ generated by $\Omega_n$ is usually called
the {\em Liouville measure}.

Now we return to the manifold $(M,\Omega)$ with a divergence-free vector field $X$ and
assume $X$ to have a smooth cross-section $N.$ This means that $X(m)\notin T_mN$ for any $m\in N$.
Suppose also for some subdomain $N_1 \subset N$ orbits of $X$ return to $N$ in
finite time. Then the Poincar\'e map $P: N_1 \to N$ is defined. This map
is a diffeomorphism from $N_1$ onto the image $P(N_1)$. As is known, Poincar\'e map $P$ on $N_1$ generated by the flow $f^t$
is defined as follows.
For a point $x\in N_1$ its image $P(x)$ is defined as $P(x)=f^{T(x)}(x),$ where $T(x)$ is the first return
time to $N$ for the orbit through $x$. Due to transversality of $N$ to the vector field at the points of $N$
and smoothness of $X$, function $T(x)$ depends smoothly on $x$, hence $P$ is a smooth map $P:N_1 \to
N$. The following assertion is used in many papers, we present the proof
for the reader convenience.
\begin{proposition}
$P$ is symplectic w.r.t. 2-form $\omega_X$ on $N$, in particular, $\omega_X$ is nondegenerate.
\end{proposition}
\proof To prove this fact, let us remark first the 2-form on $\omega_X$ as a function of two vectors
$\omega_X(\xi,\eta)$ is bilinear and skew-symmetric. This follows from the
properties of $\Omega.$ It is also nondegenerate.
Indeed, let for a fixed $\xi \in T_xN$ one has $\omega_X(\xi,\eta)=0$
for any $\eta \in T_xN.$ We need to prove that $\xi=0.$ If $\xi \ne 0,$ consider the plane $L_\xi$
in $T_xM$ spanned by two independent vectors $\xi$ and $X(x)$. If $\eta \in
L_\xi$, then $\Omega(X(x),\xi,\eta)=0$, but if $\eta \notin L_\xi$, then
vectors $X(x), \xi, \eta$ are not coplanar and $\Omega(X(x),\xi,\eta)\ne
0$ since $\Omega$ is the volume form. But intersection of two transverse planes $T_xN$ and $L_\xi$
is the straight-line
$l_\xi \subset T_xN$ spanned by $\xi$. So, $\omega_X(\xi,\eta)=0$ implies $\eta \in l_\xi$.
But $\eta$ is, by assumption, an arbitrary vector in $T_xN$, so $\xi = 0.$ Also, on
2-dimensional manifold $N$ any nondegenerate skew-symmetric
2-form is closed.

Let us shift the triple $(X(x),\xi,\eta)$ to $T_yM$ by the differential of the flow,
where $y=f^{T(x)}(x)\in N.$ We get $(Df^TX(x),Df^{T}\xi,
Df^T\eta)$, here $T(x)$ is the transition time for point $x\in N_1$. Then the equality holds
$$
\Omega(X(x),\xi,\eta)= (f^T)^*\Omega(X(x),\xi,\eta)=\Omega(Df^T(X(x)),Df^{T}(\xi),
Df^T(\eta)),
$$
since flow maps preserve $\Omega$. Denote $\omega$ the restriction of the form $\omega_X$ on $N$.

To verify the symplecticity of the Poincar\'e map $P$ we
calculate the differential $DP$ at $x$. Let us take on $N$ a smooth curve $c(s)$
through $x$, $c(0)=x,$ $c'(0)=\xi,$ and transform it by the map
$P(c(s))=f^{T(c(s))}(c(s)) \in N$. Differentiating the transformed curve and setting $s=0$
gives
$$
DP(\xi)= Df^{T(x)}(\xi) + X(P(x))dT_x(\xi).
$$
Thus the difference between $DP(\xi)$ and $Df^{T(x)}(\xi)$ is proportional to the vector
$X(P(x))$ with the proportionality factor $dT_x(\xi)$.

In addition, in accordance with the definition of the orbit, one has
$$
\frac{d}{d\tau}(f^{t+\tau}(x))|_{\tau=0}= \frac{d}{d\tau}(f^\tau\circ f^{t}(x))|_{\tau=0} = X(f^t(x)).
$$
Calculating this in the reverse order and using the group properties of the flow diffeomorphisms we come to
$$
\frac{d}{d\tau}(f^{t+\tau}(x))|_{\tau=0}= \frac{d}{d\tau}(f^t\circ f^{\tau}(x))|_{\tau=0} = Df^t(X(x)).
$$
Thus we have the identity for the flow
\begin{equation}\label{iden}
Df^t(X(x))\equiv X(f^t(x)).
\end{equation}

Now we come to the equality
\begin{multline*}
\omega_x(\xi,\eta) = \Omega(X(x),\xi,\eta)= (f^T)^*\Omega(X(x),\xi,\eta)=
\Omega(Df^{T}X,Df^{T}\xi, Df^{T}\eta)=\\ =\Omega(X(P(x)),DP(\xi)-dT_x(\xi)X(P(x)),DP(\eta)-dT_x(\eta)X(P(x)))= \\
=\Omega(X(P(x)),DP(\xi),DP(\eta))= \omega_{P(x)}(DP(\xi),DP(\eta)).
\end{multline*}
$\blacksquare$

A similar construction being applied to the case of divergence-free vector fields
on a smooth $n$-dimensional oriented manifold $(M,\Omega)$ with a nondegenerate $n$-form $\Omega$
and a smooth cross-section $N$ leads to the study of volume preserving
mappings on $N$ w.r.t. the nondegenerate $(n-1)$-form $\omega$. In the
case of nonautonomous $T$-periodic in $t$ divergence-free vector fields on an
oriented 3-dimensional $M$ this gives for its Poincar\'e map in period $T$ the
volume-preserving diffeomorphism of $M$.

These calculations show that if the divergence-free flow has a cross-section, then
the related Poincar\'e map is symplectic and all known results on such maps
are applicable (see, for instance, \cite{Meiss,KH}). In particular, if a divergence
free vector field $X_0$ has a smooth integral (integrable, see Sec.\ref{int}) and possesses
a domain in the phase manifold filled with invariant tori, then its perturbation obeys the
conclusions of KAM theory (see, for instance, \cite{BHT}). Such type results were elaborated
when studying of divergence-free vector fields near singularities \cite{Broer}. Also, all
structures related with the stochastic properties of these maps reflect
themselves in the chaotic behavior of the trajectories of liquid
particles. But, of course, it is not obligatory, when the flow has a global
cross-section. For instance, the well known ABC flow [4,5] on 3-torus most
likely does not have global cross-sections. On the other hand, majority of such
flows have periodic orbits, hence local cross-sections exist. In any case,
this is the main reason why structures characteristic for the Hamiltonian
flows are observed in liquid currents \cite{Ar,Henon,BSh,Kozlov,Lamb,MK,MW}.

Similar considerations in Lagrangian periodic incompressible liquid flows
lead to the study of iterations of volume preserving diffeomorphisms.
There is a vast literature devoted to this topic. Many similarities with
the structures observed in Hamiltonian dynamics can be found in these
investigations. But there are many differences since the property of volume preservation
gives much more freedom for the behavior of orbits \cite{LM,SM}.

\section{Divergence-free vector fields and Hamiltonian systems}

Now we suppose a divergence-free vector field $X$ on an oriented smooth 3-manifold
$(M,\Omega)$ to have a global cross-section. Due to a discussion above, this  means that for some
smooth closed 2-dimensional submanifold $N \subset M$ all orbits of $X$
intersect $N$ transversely and each orbit starting on $N$ returns on $N$ in a finite time.
Then the Poincar\'e map $P: N\to N$ is globally defined on $N$ and $M$ itself is a
bundle over $S^1$ with a leaf $N$. In other words, $M$ is diffeomorphic
to the suspension over a diffeomorphism $P: N\to N$
and a roof function $F: N\to \mathbb R$ being the return time $F(x)$
for the orbit through $x\in N$.

The construction above implies that there is a symplectic 2-form $\omega$
on $N$ such that $P$ is a symplectic diffeomorphism
w.r.t. $\omega: P^*\omega = \omega.$ In this section we want to show that there is a
4-dimensional smooth symplectic manifold $(\tilde M, \Lambda)$ and a
smooth Hamilton function $H$ on $\tilde M$ such that
the Hamiltonian vector field $X_H$ on some its level $H=c$ coincides with
the vector field $X$. Namely, the following theorem holds
\begin{theorem}
There is a smooth manifold $\tilde{M}$ of dimension four and a symplectic
2-form $\Lambda$ on $\tilde{M}$ such that
\begin{itemize}
\item $M$ is a smooth submanifold of $\tilde{M}$,
\item vector field $X$ is extended till a Hamiltonian vector field $\tilde{X}$
on $(\tilde{M},\Lambda)$ with a smooth Hamiltonian $H: \tilde{M}\to \R$ and $M$ is the level $H=c$ for some $c$.
\end{itemize}
\end{theorem}
\proof
It is sufficient to consider $M$ and the flow of the vector field $X$ as given by the suspension over a diffeomorphism
$P:N\rightarrow N$. To construct the manifold $\tilde{M}$ and the vector field $\tilde{X}$ we
increase  by one the dimension of the suspension construction.

Let us define the action $(x,r,s)\cdot n$ of the group $\mathbb{Z}$ on the manifold
$N\times\mathbb{R}^2$ setting $(P^n(x),r-n,s)$
for any $x\in N$ and all $n\in\mathbb{Z}$. This action is free and its set of orbits
$\tilde{M}=(N\times\mathbb{R}^2)/\mathbb{Z}$ is a smooth manifold
of dimension four and the quotient mapping $q:N\times\mathbb{R}^2\rightarrow\tilde{M}$
is a regular covering.

If $\mathbb{R}^2\times\mathbb{Z}\rightarrow\mathbb{R}^2$ is the auxiliary action defined as $(r,s)\cdot
n=(r-n,s)$, then the factor-manifold $C=\mathbb{R}^2/\mathbb{Z}$ is diffeomorphic to the cylinder and the quotient mapping
$q_0:\mathbb{R}^2\rightarrow C$ is the regular covering as well.

The formula $p(q(x,r,s))=q_0(r,s)$ defines a submersion $p:\tilde{M}\rightarrow C$.
Hence, a smooth bundle $\xi=(\tilde{M},p,C)$ is defined.

Let $p_0:N\times\mathbb{R}^2\rightarrow\mathbb{R}^2$ be a natural projection.
According to the construction, the following diagram is commutative
\begin{equation}\label{CD}
\begin{CD}
N\times\mathbb{R}^2 @>{p_0}>> \mathbb{R}^2\\
@VqVV @VV{q_0}V\\
\tilde{M} @>{p}>> C.
\end{CD}
\end{equation}

Consider a flow $\hat{h}^t$ on $N\times\mathbb{R}^2$ being defined by the formula $\hat{h}^t(x,r,s)=(x,r+t,s)$.
This flow generates  a vector field $\hat{X}=\partial/\partial{r}$. Since the flow $\hat{h}^t$ commutes with the action of
the group $\mathbb{Z}$, it induces the flow $\tilde{h}^t$ на $\tilde{M}$.
Denote its vector field as $\tilde{X}$.

A usual suspension over the diffeomorpism $P:N\rightarrow N$ is given, if everywhere in these constructions set
$s=0$. Then $M=q(N\times\mathbb{R}\times0)$ is a smooth 3-dimensional submanifold in $\tilde{M}$.
The flow $h^t$ corresponding to the vector field $X$ on $M$ is generated by the restriction of the flow $\hat{h}^t$
on $N\times\mathbb{R}\times{0}$. Therefore, one has $\tilde{X}|_{M}=X$, i.e. $\tilde{X}$ is the extension of
$X$ from the submanifold $M\subset\tilde{M}$ to the whole $\tilde{M}$.
Finally, the restriction of the bundle $\xi$ over $S^1=q_0(\mathbb{R}\times{0})\subset C$
gives the bundle $\xi_1=(M,p_1,S^1)$.

The standard symplectic form $\lambda_0=dr\wedge ds$ on $\mathbb{R}^2$ is invariant with respect to shifts
$(r,s)\rightarrow(r-t,s)$, $t\in\mathbb{R}$. In particular, this is true
for $t=n\in\mathbb{Z}$. Thus $\lambda_0$ generates a symplectic form $\lambda$ on $C$.
To further purposes, one needs to calculate the Lie derivative of the form $p^*\lambda$
in the direction of the vector field $\tilde{X}$.
\begin{lemma}\label{L_X_lambda}
The equality $L_{\tilde{X}}(p^*\lambda)=0$ holds.
\end{lemma}
\proof
By construction we get $\lambda_0=q_0^*\lambda$. Due to \eqref{CD} one has
\begin{equation}\label{q*p*lambda}
q^*(p^*\lambda)=p_0^*\lambda_0.
\end{equation}
The invariance of $\lambda_0$ with respect to the shifts implies that the form $p_0^*\lambda_0$
is preserved by the flow $\hat{h}^t$ on $N\times\mathbb{R}^2$. In virtue of $\tilde{h}^t\circ q=q\circ \hat{h}^t$ we have
$$
q^*(\tilde{h}^t)^*(p^*\lambda))=(\hat{h}^t)^*(q^*(p^*\lambda))=(\hat{h}^t)^*(p_0^*\lambda_0)=p_0^*\lambda_0=q^*(p^*\lambda).
$$
The mapping $q$ is a covering, therefore the equalities obtained imply that $(\tilde{h}^t)^*(p^*\lambda)=p^*\lambda$.
It remains to remember that the flow $\tilde{h}^t$ is generated by the field $\tilde{X}$.
$\blacksquare$

A natural projection $\mu_0:N\times\mathbb{R}^2\rightarrow N\times\mathbb{R}\times0$
induces the map $\mu:\tilde{M}\rightarrow M$ via the formula $\mu(q(x,r,s))=q(x,r,0)=q(\mu_0(x,r,s))$.
Here for any points $\hat{u}\in N\times\mathbb{R}^2$, $\tilde{u}=q(\hat{u})$ and $u=\mu(\tilde{u})$
we have
$$
D\mu(\tilde{X}(\tilde{u}))=D\mu(Dq(\hat{X}(\hat{u})))=Dq(D\mu_0(\hat{X}(\hat{u})))=X(u).
$$
Let us set $\tilde{\Omega}=\mu^*\Omega$ and $\tilde{\omega}=\imath_{\tilde{X}}\tilde{\Omega}$.
Then the preceding equalities imply $\tilde{\omega}=\mu^*\omega_X$.

We set $\Lambda=\tilde{\omega}+\varepsilon^{-1}p^*\lambda$ for some $\varepsilon \ne 0$.

\begin{lemma}\label{symplectic Lambda}
2-form $\Lambda$ on $\tilde{M}$ is symplectic.
\end{lemma}

\proof
First we prove that the form $\tilde{\omega}$ is nondegenerate on the leaves of the bundle $\xi=(\tilde{M},p,C)$.

Each leaf $F_1$ of the bundle $\xi_1=(M,p_1,S^1)$ is a cross-section for the flow generated by the vector field $X$.
By the proposition 1 this means that the form $\omega_X$ is nondegenerate on $F_1$.

Let $c=q_0(r,s)\in C$ and $c_1=q_0(r,0)$. We set $F=p^{-1}(c)$ and $F_1=p_1^{-1}(c_1)$.
Then one has $F=q(N\times{r}\times{s})$, $F_1=q(N\times{r}\times0)$ and
therefore the restriction $\mu|_{F}:F\rightarrow F_1$ is a diffeomorphism.
But then the form $\tilde{\omega}=\mu^*\omega_X$ is nondegenerate on $F$.

The nondegeneracy of the form $\tilde{\omega}$ on bundle leaves of $\xi$ and lemma 1
from \cite{LY} imply that $\Lambda$ is symplectic on $\tilde{M}$.
$\blacksquare$

Thus we construct a symplectic manifold $(\tilde{M},\Lambda)$ for which the first assertion of the theorem holds and
the vector field $\tilde{X}$ on $\tilde{M}$ that is the extension of the field $X$. It remains to prove that
 $\tilde{X}$ is Hamiltonian.
\begin{lemma}\label{closed i_X Lambda}
The form $\imath_{\tilde{X}}\Lambda$ is closed.
\end{lemma}
\proof
For any vector field $\tilde{Y}$ on manifold $\tilde{M}$ one has
$$
(\imath_{\tilde{X}}\Lambda)(\tilde{Y})=\Lambda(\tilde{X},\tilde{Y})=
\tilde{\omega}(\tilde{X},\tilde{Y})+\frac{1}{\varepsilon}(p^*\lambda)(\tilde{X},\tilde{Y}).
$$
But the following equalities hold
$
\tilde{\omega}(\tilde{X},\tilde{Y})=(\imath_{\tilde{X}}\tilde{\Omega})(\tilde{X},\tilde{Y})=
\tilde{\Omega}(\tilde{X},\tilde{X},\tilde{Y})=0,
$
hence
\begin{equation}\label{i_X Lambda}
\imath_{\tilde{X}}\Lambda=\frac{1}{\varepsilon}\imath_{\tilde{X}}(p^*\lambda).
\end{equation}
Also we have
$d(\imath_{\tilde{X}}(p^*\lambda))=L_{\tilde{X}}(p^*\lambda)-\imath_{\tilde{X}}(d(p^*\lambda))$.
Since $d(p^*\lambda)=p^*(d\lambda)=0$, the preceding formula and lemma \ref{L_X_lambda}
imply $d(\imath_{\tilde{X}}(p^*\lambda))=0$. In virtue of \eqref{i_X Lambda} lemma has been proved.
$\blacksquare$

\begin{lemma}\label{periods i_X Lambda}
Periods of the form $\imath_{\tilde{X}}\Lambda$ along any 1-cycle in the manifold $\tilde{M}$ are equal to zero.
\end{lemma}
\proof
The inclusion $M\rightarrow\tilde{M}$ is a homotopic equivalence. Then it is sufficient to
calculate periods along cycles lying in $M$. It is also evident that as a singular 1-cycle
one can understand a piece-wise smooth loop $\gamma:I\rightarrow M$, $I=[0,1]$, with the initial point
$\gamma(0)=\gamma(1)=q(x_0,0,0)$.

Suppose $\hat{\gamma}:I\rightarrow N\times\mathbb{R}\times0$ be a path with an initial point $\hat{\gamma}(0)=(x_0,0,0)$,
for which $q\circ\hat{\gamma}=\gamma$.
For any $t\in I$ as $N_t$ and $R_t$ we denote the leaves of trivial bundles
$N\times\mathbb{R}\times0\rightarrow\mathbb{R}$ and $N\times\mathbb{R}\times0\rightarrow N$,
through the point $\hat{\gamma}(t)$. Then we have $d\hat{\gamma}/dt=Y_2+Y_1$, where
$Y_2\in T_{\hat{\gamma}(t)}N_t$ and $Y_1\in T_{\hat{\gamma}(t)}R_t$.
Here we have $Dp_0(Y_2)=0$ and vectors $Y_1$ and $\hat{X}(\hat{\gamma}(t))$ are collinear.
This implies
$\imath_{\hat{X}}(p_0^*\lambda_0)(Y_k)=(p_0^*\lambda_0)(\hat{X}(\hat{\gamma}(t)),Y_k)=0$
for $k=1,2$. These equalities imply $\imath_{\hat{X}}(p_0^*\lambda_0)(d\hat{\gamma}/dt)=0$ and then
\begin{equation}\label{integral of i_X p_0 lambda}
\int_{\hat{\gamma}}\imath_{\hat{X}}(p_0^*\lambda_0)= \int_0^1\imath_{\hat{X}}(p_0^*\lambda_0)(d\hat{\gamma}/dt)=0.
\end{equation}

Due to \eqref{i_X Lambda} and \eqref{q*p*lambda} we come to
\begin{equation}\label{integral of i_X Lambda}
\int_{\gamma}\imath_{\tilde{X}}\Lambda=
\frac{1}{\varepsilon}\int_{\gamma}\imath_{\tilde{X}}(p^*\lambda)=
\frac{1}{\varepsilon}\int_{\hat{\gamma}}\imath_{\hat{X}}(q^*(p^*\lambda))=
\frac{1}{\varepsilon}\int_{\hat{\gamma}}\imath_{\hat{X}}(p_0^*\lambda_0).
\end{equation}
From \eqref{integral of i_X Lambda} and \eqref{integral of i_X p_0 lambda} the proof follows.
$\blacksquare$

Lemmata \ref{closed i_X Lambda} and \ref{periods i_X Lambda} imply, in accordance
with de Rahm theorem \cite[chapter 5]{Warner}, that
the form $\imath_{\tilde{X}}\Lambda$ is exact. Therefore there is a smooth function
$H:\tilde{X}\rightarrow\mathbb{R}$ for which the identity $dH=\imath_{\tilde{X}}\Lambda$ holds.
But this means that the vector field $\tilde{X}$ is globally Hamiltonian.

To complete the proof we need to show that $M$ is a level of the
Hamiltonian $H$. Fix some point $x_0\in N$ and set $\tilde{u}_0 = q(x_0,0,0)$. Then the value of Hamiltonian
$H:\tilde{M}\rightarrow \mathbb{R}$ at an arbitrary point $\tilde{u}=q(x,r,s)$
can be calculated by means of the formula
$H(\tilde{u})=\int_{\gamma}\imath_{\tilde{X}}\Lambda$,
where $\gamma:[0,1]\rightarrow\tilde{M}$ is some piece-wise smooth path with the end points
$\gamma(0)=\tilde{u}_0$ and $\gamma(1)=\tilde{u}$.

We put also $\tilde{u}_s=q(x_0,0,s)$ for all $s\in\mathbb{R}$. The integral
$\int_{\gamma}\imath_{\tilde{X}}\Lambda$ depends only on the extreme point of the path $\gamma$.
This allows one to regard
that $\gamma=\alpha_s\gamma_{rs}$ where $\alpha_s(t)=q(x_0,0,ts)$ and the path
$\gamma_{rs}:[0,1]\rightarrow\tilde{M}$ lies in the submanifold
$M_s=q(N\times\mathbb{R}\times{s})$, it joins points $\tilde{u}_s$ and $\tilde{u}$.
But for such a path $\gamma_{rs}$ the equalities hold: \eqref{integral of i_X Lambda}
and \eqref{integral of i_X p_0 lambda}. Hence we get $\int_{\gamma_{rs}}\imath_{\tilde{X}}\Lambda$=0
and
$$
H(\tilde{u})=\int_{\alpha_{s}}\imath_{\tilde{X}}\Lambda+\int_{\gamma_{rs}}\imath_{\tilde{X}}\Lambda=
\int_{\alpha_{s}}\imath_{\tilde{X}}\Lambda.
$$
Thus, the value $H(\tilde{u})$ does not depend on the number $r$. Therefore, submanifolds
$M_s\subset\tilde{M}$ are the level sets of the Hamiltonian $H$. In particular, there is a number
$c\in\mathbb{R}$ such that $M=M_0=\{\tilde{u}\in\tilde{M}|H(\tilde{u})=c\}$.
$\blacksquare$

It is evident that the extension of a DFVF given on smooth
3-dimensional oriented manifold $M$ till a Hamiltonian vector field on some
symplectic 4-dimensional $\tilde{M}$ is impossible if $X$ has singular
points (equilibria) \cite{LF,Gr}. Indeed, as is known, equilibria of a Hamiltonian
vector field coincide with critical points of the Hamiltonian. But usually
a critical level of a Hamiltonian is not a smooth manifold since it has
singularities at critical points (near them the level is not a smooth
manifold).

Nevertheless, the reduction to 2-dimensional symplectic
diffeomorphisms here is exploited also, for instance, when studying homoclinic
dynamics. This can be seen, in particular, in \cite{BSh}. In these cases
the related Poincar\'e maps are discontinuous but symplectic since they
are constructed on cross-sections to homoclinic orbits to an equilibrium
where the discontinuity takes place on the trace of stable (unstable) manifold to the
related saddle or saddle-focus equilibrium. Other equilibria are also
possible but they are degenerate (have either zeroth or pure imaginary eigenvalues)
because the condition $div X =0$ holds at the equilibrium.

There is an interesting case where a close construction is exploited in other
circumstances \cite{DM}. This concerns the case when a map $f:M \to M$ acts
on a smooth manifold with the and it is assumed in addition
that there is a symmetry for $f$, that is a smooth vector field $v$ on $M$ such that $f^*v =
v$ where $[f^*v](x) = Df^{-1}v(f(x)).$ Hence $f$ transforms $v$ to itself and then orbits
of the flow are transformed to the orbits of the same flow. Suppose that
the flow of $v$ is complete (all orbits are extendable onto whole $\R$) and has a global
cross-section $\Sigma$. In this case there is a covering space $\Sigma\times \R \to M$
generated by the Poincar\'e map on $\Sigma.$ It appears there exists a lift $F$ of $f$
to $\Sigma\times \R$ which has a skew product map form with the base
$\Sigma$. Moreover, if $f$ preserves a volume form $\Omega$ on $M$ and $v$
is divergence free, then the map on the base $\Sigma$ preserves the induced
form on $\Sigma.$ This shows connections with our results.

\section{On integrable divergence-free vector fields}\label{int}

In this section we single out a class of DFVF with a simple structure.
They can be taken as initial systems for applying perturbation methods.

Suppose on a smooth 3-dimensional oriented manifold $(M,\Omega)$ a
smooth DFVF $X$ be given. We shall call such a
vector field to be {\em integrable} if it possesses a smooth integral $F$, i.e. a smooth
function $F: M\to \R$ which satisfies the identity $dF(X)\equiv 0$.
A particular case of this situation was studied, in particular, in \cite{MW}. There is
also another case of integrability similar to that for a general 3-dimensional systems,
where by the integrability one understands the existence two integrals
independent almost everywhere. This case for a DFVF we call {\em
super-integrable}.

If a vector field $X$ is integrable, then the manifold $M$ is foliated into levels of this function $F=c.$
A first natural question in this case arises: do some restrictions exist on the topology of levels of
function $F$ and flows generated by $X$ on the invariant subset $F=c$?
Recall that usually the integrability of 3-dimensional vector fields requires to have
two independent (almost everywhere) integrals (what we called super-integrability above).
We shall show that for many goals it is enough to
have only one smooth integral to investigate the orbit structure on the
majority of levels of the integral. This also shows the close relation
with integrable Hamiltonian vector fields where the (Liouville) integrability of
a 2-degrees-of-freedom Hamiltonian vector field $X_H$ follows from the
existence of one additional smooth integral independent of $H$ almost
everywhere (due to Liouville-Arnold theorem \cite{Ar_sib}). But as we shall see,
the orbit behavior on levels $F=c$ can be a
bit more complicated than in the integrable Hamiltonian case: the flow on
such a level is not always linear, this depends on interrelations between
arithmetics of rotation numbers and smoothness of the flow as in the Kolmogorov theorem
\cite{Kolmogorov}. The details on the orbit structure of the integrable
Hamiltonian systems with two degrees of freedom see in \cite{LU,BF}.

Let us consider the case when $dF \ne 0$ on a level $F=c$ (more precisely, its connected
component) assuming this component be compact. We call such a level nonsingular.
Thus a smooth closed connected 2-dimensional submanifold
$\Sigma = \{F=c\}$ is an invariant set w.r.t. the flow generated by DFVF $X$.
\begin{proposition}
A closed nonsingular invariant 2-dimensional level $\Sigma$ is orientable.
If $X$ has not zeroes on $\Sigma,$ then it is an invariant torus.
Moreover, all close to $\Sigma$ levels are also smooth invariant tori with flows without zeroes.
\end{proposition}
\proof Consider a tubular neighborhood $B$ of $\Sigma$. Since $dF \ne 0$ on
it, $\Sigma$ separates $B$ into parts with different signs of the function
$F-c$. On the other hand, if $\Sigma$ would be nonorientable, then it does
not separate $B$ \cite{Nonor}. Thus, $\Sigma$ is an orientable closed
invariant 2-dimensional manifold. If $X$ has not zeroes on $\Sigma$, then
it cannot be the Klein bottle and it is a torus. In this case, all close
levels of $F$ are diffeomorphic to $\Sigma$ by the Morse theory
\cite{Milnor} and carry flows without zeroes, that is, they are smooth
invariant tori.
$\blacksquare$

Now we consider the level $\Sigma = \{F=c_0\}$ without equilibria on
$\Sigma.$ In principle, a flow without equilibria on a smooth torus can have different structures.
If the flow has a global cross-section, then its structure depends mainly
on its Poicar\'e rotation number and the smoothness of the flow
\cite{Hartman}. But the flow on a torus without a global cross-section may have a Reeb
component \cite{HF}. We want to show that the divergence-free property
implies the strong restrictions on the flow behavior. Namely, we shall
prove the existence of smooth invariant measure for the flow.

In particular, 2-form $\omega_X(\xi,\eta)$ is completely degenerate on $\Sigma$.
Indeed, if vector $\xi \in T_x\Sigma$ is not collinear to $X(x)$, then the plane $L_\xi$
in $T_xM$ spanned by two independent vectors $\xi$ and $X(x)$ coincides with the tangent plane
$T_x\Sigma$. Thus for any vector $\eta \in T_x\Sigma$ one has
$\omega_X(\xi,\eta)=0$. This is an analog for $\Sigma$ to be a Lagrangian torus in the symplectic setting.

So, in order to get an invariant measure on $\Sigma$, we need to go in another way. The natural
way is the following. Take a thin layer $F = c,$ $|c-c_0|< \varepsilon,$
where $\varepsilon$ is small enough and positive. Then these levels of $F$
are also smooth tori without equilibria of the vector fields. The flow
in this layer preserves the volume. Let us introduce some smooth Riemannian metrics
in this layer. Since $\Sigma$ is two-sidedly imbedded, then a smooth field of
normal vectors on $\Sigma$ can be found. Choose such a field and denote $n(x)\in T_xM$ its
normal vector at the point $x\in \Sigma.$ Then a 2-form $\omega_n$ on $\Sigma$, $\omega_n(\cdot,\cdot) =
\Omega(n(x),\cdot,\cdot)$ is defined.

\begin{theorem}
A vector field of normals can be chosen in such a way that: i) 2-form  $\omega_n$ on $\Sigma$ is nondegenerate;
ii) the restriction of $X$ on the level $F=c_0$ defines the flow $\varphi^t_{\Sigma}$
that preserves the form $\omega_n$, $(\varphi^t_{\Sigma})^*\omega_n = \omega_n.$
\end{theorem}

\proof
Choose $\varepsilon\in\mathbb{R}$, $\varepsilon>0,$ small enough such that on submanifold
$M_{\varepsilon}=\{x\in M||F(x)-c_0|<\varepsilon\}$ of $M$ the inequalities $dF\ne 0$ and $X\ne 0$ stay valid.

Let $g$ be some smooth Riemannian metrics on $M_{\varepsilon}$. Then a vector field $\nabla{F}$ on $M_{\varepsilon}$
is correctly defined by the relation $\imath_{\nabla{F}}g= g(\nabla F,\cdot)=dF(\cdot)$.
By construction, one gets $\nabla{F}\ne 0$ on $M_{\varepsilon}$. Therefore, the function
$\rho=1/g(\nabla{F},\nabla{F})$ and a vector field $n=\rho\nabla{F}$ are defined and differentiable. Observe that
\begin{equation}\label{dF(n)}
dF(n)=(\imath_{\nabla{F}}g)(n)=g(\nabla{F},n)\equiv 1.
\end{equation}
Let us set as above $\omega_n=(\imath_n\Omega)|_{\Sigma}$.

The vector field $X$ is tangent to $\Sigma$, therefore we get
\begin{equation}\label{L_X omega_n}
L_X\omega_n=(L_X(\imath_n\Omega))|_{\Sigma}.
\end{equation}
By definition $\imath_n\Omega=c(n\otimes\Omega)$, where $c$ denote the contraction of the tensor
$n\otimes\Omega$. In accordance to proposition 3.2 in
\cite{KN} we have
$$
L_X(c(n\otimes\Omega))=c(L_X(n\otimes\Omega))=c(L_Xn\otimes\Omega)+c(n\otimes L_X\Omega),
$$
where $L_Xn=[X,n]$ is the Lie bracket of vector fields. As a result we come to the relation
\begin{equation}\label{L_X i_n Omega}
L_X(\imath_n\Omega)=\imath_{[X,n]}\Omega+\imath_n(L_X\Omega).
\end{equation}

\begin{lemma}\label{[X,n]}
The vector field $[X,n]$ is tangent to $\Sigma$.
\end{lemma}
\proof The vector $n(x)$ is orthogonal to the surface $\Sigma$, thus the representation $[X,n](x)=hn(x)+Z$ holds,
where $h\in\mathbb{R}$ and $Z\in T_x\Sigma$. But then, due to
\eqref{dF(n)} we get
$$
[X,n](x)F=hn(x)F+ZF=hdF(n(x))+dF(Z)=h.
$$

From the other hand, one has $nF=dF(n)\equiv1$ and $XF=dF(X)\equiv 0$.
Therefore, the equality holds
$$
[X,n]F=X(nF)-n(XF)\equiv 0.
$$
The equalities obtained lead to $h=0$ and hence $[X,n](x)=Z\in T_x\Sigma$.
$\blacksquare$

Returning to the proof of the theorem, consider arbitrary $x\in\Sigma$ and $Y,Z\in T_x\Sigma$.
It follows from the equality $\dim{\Sigma}=2$ and lemma \ref{[X,n]} that three vectors $Y,Z$ and $[X,n](x)$
are coplanar. Therefore we get $(\imath_{[X,n]}\Omega)(Y,Z)=\Omega([X,n](x),Y,Z)=0$
and hence
\begin{equation}\label{i_[X,n] Omega}
(\imath_{[X,n]}\Omega)|_{\Sigma}=0.
\end{equation}

Since $L_X\Omega=0$, then from \eqref{L_X i_n Omega} and \eqref{i_[X,n]
Omega} the equality
$(L_X(\imath_n\Omega))|_{\Sigma}=0$ follows and hence $L_X\omega_n=0$, in virtue of \eqref{L_X omega_n}.

The orthogonality of $n$ to the surface $\Sigma$ implies also that the form $\omega_n$
is nondegenerate and consequently is the area form on $\Sigma$.$\blacksquare$

Thus, for any nondegenerate level $\Sigma_c = \{F=c\}$ we have the vector field
$X_c$ being the restriction of $X$ onto the invariant submanifold
$\Sigma_c$. This vector field is nonsingular and preserves the area form
$\omega_n$. Let us choose some angle variables $(\varphi,\psi)$ on the
torus $\Sigma.$ Then 2-form $\omega_n$ takes the form $a(\varphi,\psi)d\varphi\wedge
d\psi$ with the smooth doubly periodic positive $a$ and the vector field has the form
$$
\dot\varphi = A(\varphi,\psi),\;\dot\psi = B(\varphi,\psi),
$$
where $A^2 + B^2 \ne 0$ and both smooth functions $A,B$ are doubly periodic.
The measure preservation means here
$$
\frac{\partial}{\partial \varphi}(aA)+\frac{\partial}{\partial \psi}(aB)=0.
$$
Denote
$$
\lambda_1 = \int\limits_{\Sigma}Aad\varphi\wedge d\psi,\;\lambda_2 = \int\limits_{\Sigma}Bad\varphi\wedge d\psi.
$$
The number $\lambda = \lambda_1/\lambda_2$ called the Poincar\'e rotation number plays the main role in the orbit
dynamics on the torus $\Sigma$. As is known, if $\lambda$ is
rational or one of $\lambda_i$ is equal to zero, then all orbits of the
flow are periodic (this is because of the existence of a smooth invariant
measure). But if $\lambda$ is irrational and the flow is of smoothness
$C^2$ then all orbits on the torus are transitive. More subtle ergodic properties of the flow
depend on the interrelations between the arithmetic type of $\lambda$
and a smoothness of functions $A,B$ \cite{Kolmogorov}. For  instance, for the case of $C^5$-smooth
r.h.s. the flow can have a continuous spectrum \cite{Katok} (also see details in
\cite{KSF}).

One can think that the integrability of $X$ imply the flow structure like
in the integrable Hamiltonian system. In fact,
the existence of an integral in $M$ does not imply that this additional
integral can be extended onto the 4-dimensional symplectic manifold
constructed above. This is not the case even for the case of a Hamiltonian
system with two degrees of freedom. Such a system can be nonintegrable in the
whole phase space but be integrable on some separate level of Hamiltonian.
To present such an example, let us consider some smooth symplectic
four-dimensional manifold $(M,\Lambda)$ with a symplectic form $\Lambda$
and a smooth Hamiltonian $H_0$. We assume the related Hamiltonian vector field
$X_{H_0}$ be integrable, that is there is an additional smooth integral
$K$, $\{H_0,K\}\equiv 0.$ Consider a perturbation of this vector field $H=H_0 + \eps
H_1.$ Let us fix $c$. One can choose the function $H_1$ in the form $H_1 = (H_0 - c)F$ such
that on the level $H_0 = c$ the integrable system has some integrable
structure, and a function $F$ can be taken arbitrarily. Let $J_x: T_xM \to
T^*_xM$ be the isomorphism between 1-forms and vector fields on $M$ defined by the symplectic form
$\Lambda.$ Then $JdH$ is the Hamiltonian vector field generated by function
$H$. Thus we get the Hamiltonian vector field
$$
J(dH_0 + \eps FdH_0 + \eps(H_0 - c)dF).
$$
On the set $H = c$ we have $(H_0 - c)(1+\eps F)=0$, thus it coincides with the level of
the function $H_0$ and therefore is the invariant submanifold where the
dynamics is integrable. Indeed, $X_H= J(1+\eps F)dH_0$ on this level, hence the vector field is
obtained by the change time from the integrable vector field $X_{H_0}$
on the level $H_0 = c$. It is evident that
function $F$ can be chosen in such a way that the complete dynamics would
be nonintegrable. For instance, it can be reached, if on the level $H_0 = c$
the orbit structure will contain a saddle periodic orbit with its merged stable
and unstable manifolds. Function $F$ can be chosen in such a way that on close
levels stable and unstable manifolds of saddle periodic orbits existing when changing $H=c$ would be split transversely.
This shows that the integral $K$ generally cannot be extended onto the whole phase space.

\begin{Rm}
All constructions above where they appeal to the Riemannian
metrics use a possibility the construct a smooth Riemannian metrics on a
smooth manifold $M$, applying, for instance, a partition of unity \cite{KN}.
If the manifold $M$ is real analytic this construction does not work
and one needs to use other tools. For example, to find a Riemannian
metrics in this case one can exploit the Morrey-Grauert theorem about
an analytic embedding of $M$ into the Euclidean space (see, for instance,
\cite{Sh}) and then to restrict the Euclidean metrics to this embedded
manifold. Thus all constructions can be done analytic as well. In particular,
this concerns the case of an integrable DFVF and the flow on its
nonsingular level without zeroes. The flow in this case is real analytic
and all conclusions of the Kolmogorov theorem hold.
\end{Rm}

\subsection{On global structure of integrable divergence free vector fields}

For the case of integrable DFVFs one can develop a global theory of such vector fields
similar to the case of integrable Hamiltonian vector fields on smooth symplectic four-dimensional
manifolds \cite{BF,LU}. Recall that if $H$ be a smooth Hamilton function and $F$
its additional integral, one can restrict the Hamiltonian vector field $X_H$ on some
nondegenerate level $V_c = \{H=c\}$ where $dH\ne 0$.
Suppose the restriction $F_c$ of this additional integral be a smooth Bottian function \cite{BF}, that is it has
finitely many critical values and the related critical sets of $F_c$ are organized into the
finitely many critical closed smooth curves such that the restriction of
$F_c$ on the transverse disk to the critical curve $l$ generates a Morse function
with a nondegenerate critical point at the trace of $l$ on the disk. For the restriction of the flow
$X_H$ on $V_c$ these critical curves are usually closed periodic orbits for $X_H$ and for
Bottian integral they can be only of two types: elliptic and hyperbolic
ones. If the level $V_c$ is a closed manifold (compact without a boundary), then almost all levels of $F$
are invariant Lagrangian tori (by the Liouville-Arnold
theorem). This foliation can be described via the invariant introduced by
Fomenko (see details in \cite{BF}).

A similar theory can be also elaborated for the case of integrable DFVFs.
Suppose a smooth closed oriented 3-dimensional manifold $(M,\Omega)$ be given and $X$ be a DFVF
which we assume to be integrable and without equilibria (since $M$ is
oriented and closed, its Euler characteristic is zero).
Let $F$ be the related smooth integral, i.e. $dF(X)\equiv 0,$ such
that it has a finite number of critical levels. We assume that each critical set of $F$
consists of finitely many smooth disjoint closed curves $l_1,\ldots,l_m$ such that on a cross-section to such
curve $l_i$ the restriction of $F$ is a smooth function that has a Morse critical point at the trace of
$l_i$. A simplest example of such situation is an integrable
Hamiltonian vector field on a symplectic manifold restricted on the
non-degenerate level of its additional integral when this restriction is a
Bottian function. There are many such examples in mechanics (see, for instance, \cite{BF,CB}.

Suppose now $\gamma$ be a periodic orbit of $X$ and $\{F=F(\gamma)=f\}$
be a connected set containing $\gamma.$ If this level is not critical,
then its component containing $\gamma$ is an invariant smooth torus with a rational rotation number
since it contains a periodic orbit. Due
to the existence of smooth invariant measure all orbits of $X$ on this
torus are periodic and $\gamma$ cannot be isolated. Thus, if $\gamma$ is an isolated periodic orbit of $X$,
then the connected set of $F=f$ either coincides with $\gamma$ and $\gamma$ is an elliptic periodic orbit for $X$
or $\gamma$ is a saddle periodic orbit and the connected set of $F=f$ coincides with the merged stable
and unstable manifolds of $\gamma$ (recall that $F$ is preserved along any orbit of $X$).

In the first case $\gamma$ is a closed curve of maximums (minimums) for
$F$ and for the second case $\gamma$ is the curve of saddle critical
points for $F.$ This is easily derived from the Bott property like in
\cite{Fom}.

By analogy with the Hamiltonian case we shall call such an
integral to be Bottian. As we know, all noncritical levels of $F$ are disjoint invariant tori.
We assume, in addition, that almost all of these tori carry flows with irrational rotation numbers.
This assumption prevents the super-integrability --
the existence of the second independent almost everywhere integral for $X$.

First of all, we observe that an orbit of $X$
through a critical point of $F$, since this orbit is not an equilibrium of $X$,
consists of critical points of $F$. Indeed, by definition, one has $(D\Phi^t)^*(dF)=
dF(D\Phi^t).$ If $dF(\xi)=0$ for any $\xi \in T_xM$ (the point $x$ is
critical), then $dF(\eta)=0$ for any $\eta = D\Phi^t(\xi)$ at the point
$y=\Phi^t(x).$ Since $\Phi^t$ is a diffeomorphism, then $D\Phi^t$ is the
isomorphism of of tangent spaces. Hence, each connected critical set of $F$ being a smooth closed curve is a
periodic orbit of $X$. Denote such orbit $\gamma$. Choose a local
transverse disk $D$ to $\gamma$ at some its point $m\in \gamma$ and consider a smooth
function $f$ being the restriction of $F$ on $D$. Point $m$ is a nondegenerate (Morse) critical
point for $f$ (the Bott property). Hence, $m$ can be of two types, a center
or a saddle. Since $F$ is the integral, then for the case of center $\gamma$ is
enclosed by the family of invariant tori for $X$. If $m$ is a saddle, then
through $m$ two smooth segments pass composing together the intersection of $D$ with the local level
of the set $F=F(m).$ Generically, for the saddle case the related periodic orbit of $X$ through
$m$ is a saddle periodic orbit, hence one segment generates a local stable
manifold of the orbit $\gamma$ and another segment does a local unstable
manifold of $\gamma$.

As is known, in the whole phase space a saddle periodic
orbit can be of two types: orientable and nonorientable. For a periodic
orbit of $X$ the flow generates on the cross-section $D$ the Poincar\'e map $P: D \to D$ being a
symplectic mapping w.r.t. the restriction of $\Omega$ on $D$ (see above).
This symplectic map has a fixed point at $m$. The multipliers of $DF|_m$
can be generically either two complex conjugate numbers on the unit circle
(elliptic fixed point) or two different nonzero real numbers $\mu,
\mu^{-1}$ (saddle fixed point).
An orientable periodic orbit corresponds to the positive multipliers $\mu,
\mu^{-1}$ of $DF|_m$ and nonorientable periodic orbit corresponds to the negative multipliers of
$DF|_m$. Local stable and unstable 2-dimensional manifolds $W^s(\gamma),$ $W^u(\gamma)$ for the orientable periodic
orbits are both cylinders and for the nonorientable periodic orbit they are both M\"obius
strips. Their continuation by the flow of $X$ gives global stable and
unstable manifolds which we also denote as $W^s(\gamma),$ $W^u(\gamma).$
Consider the $\omega$-limit set of an orbit from $W^u(\gamma).$
Since $M$ is closed, this set is not empty, invariant and closed.

To ease the exposition, we assume that all critical sets of $F$ are either
elliptic periodic orbits (critical sets of maximums or minimums) or saddle
periodic orbits for $X$. Then for the case of Bottian integral all orbits
of $X$ and their limit sets in $M$ are known. They are
either closed periodic orbits composed by critical sets of $F$ or they
belong to stable (unstable) manifolds of saddle periodic orbits (saddle critical
curves) or they belong to invariant 2-tori. In the latter case their
$\omega$- ($\alpha$-) limit sets belong to the same torus as well.
In the case under consideration, all limit set for an orbit in
$W^u(\gamma)$ is one of saddle periodic orbits. More precisely, the
following assertion is valid
\begin{lemma}
Consider a critical level $\{F=F(\gamma)\}$ being a compact set in $M$
and let $\gamma_1,\ldots, \gamma_s$ be all saddle periodic orbits of $X$ in this set.
Then any connected component of the set $\{F=F(\gamma)\}\setminus
\{\gamma_1,\ldots,\gamma_s\}$ consists of homo- or heteroclinic orbits of
one or two different saddle periodic orbits from the collection $\{\gamma_1,\ldots,
\gamma_s\}$.
\end{lemma}
This lemma means that all stable and unstable manifolds of the same or
different saddle periodic orbits coalesce.

Now one may consider, as in \cite{Fom}, the evolution of tori and their reconstructions when
a natural parameter of the system, the value of integral $F$, varies. Since
$M$ is closed, $F$ takes its maximal $f_+$ and minimal $f_-$ critical values.
For a Bottian integral $F$ we have the related critical sets -- the collection of
elliptic periodic orbits. Each maximal (and minimal) periodic orbit $\gamma^+_i$ ($\gamma^-_j$) gives rise to the
one-parameter family of tori starting from $\gamma^+_i$ (if $f_+$ corresponds
to several maximal critical curves, then we get the related number of families of
tori). Tori from different families (or even from the same family) can
collide when they approach to the critical levels of saddle critical
curves, reconstruct and continue forming other families. The complete
reconstructions can be described by the invariant like "moleculas" (see details in
\cite{BF}).

\begin{Rm} If a DFVF in question has equilibria, then its integrability
with a smooth integral is in question. For instance, if this integral would
have nondegenerate (of the Morse type) critical points, then the local structure
of levels for this integral is determined by the Morse lemma
\cite{Milnor}. The singular level containing the critical point is
nonsmooth at the critical point (it has a cone-like singularity). But this
level is an invariant set of the flow and usually this point is simultaneously
a singular point of the vector field. Therefore, either this singular
point  has to be degenerate or the integral should have more degenerate
critical point. Examples of (super-)integrable vector field in \cite{NVV}
show this. Another possibility is a nondegenerate (hyperbolic) singular
point of a DFVF and its smooth integral with lines of critical points
being not Bottian.
\end{Rm}

\section{Acknowledgement}

The idea of this paper appeared when one of the authors (L.L.) was under preparation of his
lectures at the winter school ``Dynamical Systems and Fluid Motions'' (the University of Bremen, March 27-31,
2017, see http://wis-fluids.math.uni-bremen.de/Abstracts.html). L.L. is thankful to organizers,
especially to J. Rademacher and I. Ovsyannikov for the invitation and a care.

The authors thank F. Laudenbach for the useful discussion and
explanations as well as J. Meiss for remarks and criticism having
been allowed to improve the text. We are also thankful to D.V. Alekseevsky
for the explanations related with real analytic situation.

A financial support from the Russian Science Foundation (grant 14-41-00044) is
acknowledged, as well as from Russian Ministry of Science and Education
(project 1.3287.2017, target part) and the Russian Foundation of Basic Research
(grant 16-01-00312а).

\end{document}